\documentclass[a4paper,11pt,twoside,dvips]{article}
\usepackage{latexsym}
\usepackage{amsmath}
\usepackage{amssymb}
\usepackage{textcomp} % Euro Zeichen
\usepackage[ansinew]{inputenc} % dito
\usepackage[dvips]{graphicx} % Zum Einfgen von *.eps Grafiken (anschl. *.dvi Format, nicht *.pdf)
\usepackage{tabularx}  % Spezielle Tabellen: begin{tabularx}{L%G�%@ge}{Ausrichtung (X)}
\usepackage{amsfonts}
\usepackage{float}
\usepackage{amssymb}
\usepackage{amsthm}
\usepackage{verbatim}
\usepackage{graphicx}
\usepackage{bbm}
\usepackage[T1]{fontenc}
\usepackage{fancyhdr}
\newtheorem{thm}{Theorem}[section]

\newtheorem{lemma}[thm]{Lemma}

\newtheorem{rem}[thm]{Remark}

\renewcommand{\cal}{\mathcal}
\renewcommand{\rm}{\mathrm} % nicht kursiv
\newcommand{\bb}{\mathbb} % Doppelter Strich (reelle Zahlen etc.)
\newcommand{\df}{\stackrel{\mathrm{def}}{=}}
\newcommand{\eps}{\varepsilon}

\numberwithin{equation}{section}

\begin{document}
\title{\textbf{Special Examples of Diffusions in Random Environment}}
\author{Ivan del Tenno\\ Department of Mathematics\\ ETH Zurich\\ CH-8092 Zurich\\ e-mail: deltenno@math.ethz.ch\\
Phone: +41 44 632 3443\\[11pt]
May, 2008}
\date{\empty}
\maketitle
\begin{abstract}
In this note we present some examples of diffusions in random environment whose asymptotic behavior is rather
surprising. We construct a family of diffusions that are small perturbations of Brownian motion with non-vanishing expected
local drift under the static measure of the environment but where the ballistic behavior is lost. As slight modifications 
of this collection of diffusions we also 
provide examples with ballistic behavior where the non-vanishing limiting velocity
points to a direction opposite to the expected local drift under the static
measure.
\end{abstract}
{\textbf{Keywords:}} Diffusions in random environment, invariant ergodic measure, environment viewed from the particle\\
{\textbf{MSC:}} 82D30; 60K37

\section{Introduction}\label{a}

In this note we construct examples of diffusions in random environment which are small 
perturbations of Brownian motion and behave diffusively although their expected local drift does not vanish. 
We also provide examples where the diffusion has non-degenerate asymptotic velocity but the average local drift vanishes, as 
well as examples where the asymptotic velocity points in the opposite direction to the average local drift. 
These examples show that the naive guess that the average local drift controls whether the diffusion is 
ballistic (i.e. has non-degenerate velocity) or not, is wrong, even for small perturbations of Brownian motion. 
The crucial tool in the construction of the above examples is the so-called method of 
{\textit{the environment viewed from the particle}}, see below (\ref{62}) for more comments 
and references on this technique. 
%However this method does not seem to work in a number of situations. 
Let us mention that there has recently been an intense activity around the investigation of diffusions as well as 
random walks in random environment, notably in the study of ballistic behavior, see
\cite{goer2}, \cite{goer}, \cite{schmitz1}, \cite{schmitz2}, \cite{shen}, \cite{sznit1}, \cite{sznit4}, 
\cite{sznit2}, \cite{sznitzer}. As for diffusive behavior, some progress has also been made, 
see \cite{boltsznitzei}, \cite{del}, 
\cite{sznitzei}. For an overview of results and useful techniques concerning this area of research we also refer to
\cite{olla2}, \cite{sznit3}, \cite{sznit5}, \cite{zei1}, \cite{zei2}.\\

We now describe the model. We consider dimensions $d\geq 2$ and the random environment is described by a
probability space $(\Omega,\cal{A},\bb{P}).$ We write $E^{\bb{P}}$ for the expectation with respect to
the measure $\bb{P}.$ We assume the existence of a group $\{\tau_x\, :\, x\in\bb{R}^d\}$ of $\bb{P}$-preserving
transformations which act ergodically on $\Omega$ and are jointly measurable in $x\in\bb{R}^d,\ \omega\in
\Omega.$ Note that in the sequel each measurable function of the form $f(x,\omega),\ x\in\bb{R}^d,
\ \omega\in\Omega,$
is supposed to be generated from a measurable map $\tilde{f}$ on $\Omega$ by the action 
$\{\tau_x\, :\, x\in\bb{R}^d\},$ i.e.
\begin{equation}\label{1}
f(x,\omega)\df \tilde{f}(\tau_x(\omega)),
\end{equation}and hence is a stationary random field. For all such functions $f$ and all 
Borel subsets $F\subseteq \bb{R}^d$ we define the $\sigma$-algebra 
\begin{equation}\label{300}
\cal{H}^{f}_F\df\sigma\left(
f(x,\omega)\, :\, x\in F\right).
\end{equation}
The diffusions under consideration in this work all have 
the identity as diffusion matrix and the drift is a stationary uniformly bounded function 
$b(x,\omega),\ x\in\bb{R}^d,\ \omega\in\Omega,$ such that there exists a constant
$\kappa>0$ such that for all $x,y\in\bb{R}^d,\ \omega\in\Omega,$ the following Lipschitz 
condition is satisfied:\begin{equation}\label{2}
|b(x,\omega)-b(y,\omega)|\leq \kappa|x-y|,
\end{equation} where $|\cdot|$ denotes the Euclidean norm in $\bb{R}^d.$ Furthermore, the coefficient
$b$ satisfies the finite range dependence condition with range $R>0$
that is for Borel subsets $A$ and $B$ of $\bb{R}^d,$ 
\begin{equation}\label{4}
\cal{H}_A^{b}\mbox{ and }\cal{H}_B^{b}\mbox{ are }\bb{P}\mbox{-independent whenever }d(A,B)>R,
\end{equation}where $d(A,B)\df\inf\left\{|x-y|\, :\, x\in A,y\in B\right\}.$ We write $(X_t)_{t\geq 0}$
for the canonical process on $C(\bb{R}_+,\bb{R}^d),$ the space of continuous $\bb{R}^d$-valued functions
on $\bb{R}_+.$ Due to boundedness and regularity of $b$ (see (\ref{2})), for any $x\in\bb{R}^d,\ \omega\in\Omega,$
the martingale problem attached to \begin{equation}\label{5}
\cal{L}^{\omega}=\frac{1}{2}\Delta+b(\cdot,\omega)\cdot\nabla
\end{equation}and starting in $x$ at time 0 is well-posed, see Corollary 5.4.29 in \cite{ks} for the existence of 
a unique solution of the above martingale problem. The law $P_{x,\omega}$ on
$C(\bb{R}_+,\bb{R}^d)$ denotes its unique solution and describes the diffusion in the environment
$\omega$ and starting from $x$. $P_{x,\omega}$ is usually called the {\textit{quenched}} law and we write $E_{x,\omega}$
for the corresponding expectation. For the study of the asymptotics of $(X_{t})_{t\geq 0},$ it is convenient to
introduce also the {\textit{annealed}} law \begin{equation}\label{6}
P_x\df\bb{P}\times P_{x,\omega},\quad x\in\bb{R}^d.
\end{equation}We denote with $E_x$ the corresponding expectation.\\

Let us now explain more precisely the purpose of this work. We are going to construct in dimension $d\geq 2$ a family of
small perturbations of Brownian motion attached to a second order elliptic operator of the 
form (\ref{5})
which contains examples of diffusions with arbitrarily small drifts such that the expected local drift 
under the static measure does not vanish but the ballistic behavior is lost. More precisely, for all 
$\eps>0$ small enough we find examples with $|b(x,\omega)|\leq \eps$ for all $x\in\bb{R}^d,\omega\in\Omega,$ 
such that\begin{equation}\label{60}
E^{\bb{P}}[b(0,\omega)]\ne 0\quad\mbox{ but }\quad P_{0}\mbox{-a.s., }\lim_{t\to\infty}\frac{X_t}{t}=0,
\end{equation} see Theorem \ref{p1}.
As a slight modification of this class of diffusions, we will also provide examples with arbitrarily small drifts that 
exhibit ballistic 
behavior when the expected local drift vanishes, that is\begin{equation}\label{61}
E^{\bb{P}}[b(0,\omega)]= 0\quad\mbox{ but }\quad P_{0}\mbox{-a.s.,
}\lim_{t\to\infty}\frac{X_t}{t}=v\ne 0,
\end{equation} see Theorem \ref{p3}, with a deterministic velocity $v\ne 0,$ or even has the opposite direction to the
non-vanishing limiting velocity which means that there exists a positive constant $\gamma>0$
such that\begin{equation}\label{62}
P_0\mbox{-a.s., }\lim_{t\to\infty}\frac{X_t}{t}=-\gamma E^{\bb{P}}[b(0,\omega)]\ne 0,
\end{equation}see Theorem \ref{p2}.
For the construction of examples with similar behavior in the case of random walks in random 
environment, the authors of \cite{boltsznitzei} exploit the presence of so-called 
cut times to derive a law of large numbers, see also del Tenno \cite{del} for a
similar technique in the continuous space-time setting. In this work
the main strategy to provide examples of the nature described in
(\ref{60})-(\ref{62}) is the method of {\textit{the environment
viewed from the particle.}} 
For further successful applications of this method see for instance \cite{kipva}, \cite{komol}, \cite{koz}, 
\cite{o}, \cite{pava2}, \cite{ras}, \cite{remi2}. 
%For a successful application of this method see for instance \cite{kipva}, 
%\cite{komol}, \cite{koz}, \cite{o}, \cite{pava1}, \cite{pava2}, \cite{remi1}, \cite{remi2}. 
This technique relies on the existence of an invariant ergodic measure for the process of the environment 
viewed from the particle, as defined in (\ref{11}), which is absolutely continuous with respect 
to the static measure of the environment and produces a law of large numbers 
with explicit formula for the limiting velocity. Unfortunately, for general diffusions in random environment,
the problem of finding such a measure seems to be intractable. For our purpose however,
 we can restrict ourselves to drifts with a special structure for which the invariant
 measure is known, see (\ref{01}) and Theorem \ref{t01} in the Appendix.\\
The construction of the above mentioned examples using the technique of {\textit{the environment viewed from 
the particle}} shows us that the limiting velocity of the particle is
governed by the environment viewed from the particle and not by the static environment which may be different.\\

Let us explain how this work is organized. In Section \ref{b} we construct a family
of small perturbations of Brownian motion which contains examples of each type of behavior described above in (\ref{60})-(\ref{62}).
In Section \ref{c} we then provide these examples. Section \ref{e} contains auxiliary results and 
is dedicated to the proof of Lemma \ref{l3} stated in 
Section \ref{b} and a concrete example of possible functions which are involved in the construction of the family of 
drifts in Section \ref{b}. Finally, in the Appendix, we will present more details
on the existence of an invariant ergodic measure for the process of the environment viewed from the
particle for a special class of diffusions.\\

%{\textbf{Convention on constants:}} Unless otherwise stated all the generic constants in this article depend on 
%the dimension $d$, the uniform bound on the drift $\sup_{x\in\bb{R}^d,\omega\in\Omega}|b(x,\omega)|<\infty$ and 
%the Lipschitz constant of the drift $b$, see (\ref{2}), of the corresponding diffusion under 
%consideration.\\

\section{Main construction of the drifts}\label{b}

In this section we are going to construct a family of small drifts 
%that are uniformly bounded and satisfy the
%regularity condition (\ref{2}) and the finite range dependence condition with range 
%$R>0,$ see (\ref{4}), 
such that the class 
of diffusions generated via the operators given in (\ref{5}) for this family 
of drifts contains examples of diffusions with arbitrarily small drift which behave as described in the introduction 
(see (\ref{60})-(\ref{62})).\\

Let us first introduce some notation. For a Borel subset $A\subseteq \bb{R}^d$ and a real number $r>0$
we define the open $r$-neighborhood of $A$ as\begin{equation}\label{76}
A^r\df\left\{ x\in\bb{R}^d\,\Big|\,d(x,A)< r\right\}, 
\end{equation} with the distance function as defined below (\ref{4}). For the $j$-th, $j\geq 0,$ partial derivative with 
respect to $x\in\bb{R}^d$ of a measurable function $f(x,\omega),\ x\in\bb{R}^d,\ \omega\in\Omega,$ 
we will write\begin{equation}\label{19}
\partial_i^jf(x,\omega)\df \frac{\partial^j}{\partial x_i^j}f(x,\omega),\quad i=1,\ldots,d.
\end{equation}

For simplicity we will only write $\partial_i$ for $\partial_i^1,\ 1,\ldots,d.$
Furthermore, for integer $k\geq 0,$ we denote with $C^k(\bb{R}^d)$ the space of all 
$k$-times continuously differentiable real-valued functions on $\bb{R}^d$, whereas 
$\mathrm{Lip}^k_{m}(\bb{R}^d\times\Omega)$ for $m\geq 0$ stands for the space of all 
measurable functions $f:\bb{R}^d\times\Omega\longrightarrow \bb{R}$ with 
$f(\cdot,\omega)\in C^k(\bb{R}^d)$ for each $\omega\in\Omega$ and such that for all 
$x,y\in\bb{R}^d,\ \omega\in\Omega$ and $i=1,\ldots,d,\ j=0,1,\ldots,k,$\begin{equation}\label{20}
|\partial_i^jf(x,\omega)|\leq m,\quad |\partial_i^jf(x,\omega)-\partial_i^jf(y,\omega)|
\leq m|x-y|.
\end{equation}Note that no such control on mixed derivatives is needed in the sequel.\\

%For simplicity we will only write $\partial_i$ for $\partial_i^1,\ 1,\ldots,d.$ For a 
%multi-index $J=(j_1,\ldots,j_d)$ with $j_i\in\{0,1,\ldots\}$ we define its length as $|J|=
%j_1+\cdots+j_d$ and denote the corresponding mixed derivatives by\begin{equation}\label{0211}
%\partial^Jf(x,\omega)\df\partial^{j_1}_1\cdots\partial_d^{j_d}f(x,\omega).
%\end{equation}
%Furthermore, for integer $k\geq 0,$ we denote with $C^k(\bb{R}^d)$ the space of all 
%$k$-times continuously differentiable real-valued functions on $\bb{R}^d$, whereas 
%$\mathrm{Lip}^k_{m}(\bb{R}^d\times\Omega)$ for $m\geq 0$ stands for the space of all 
%measurable functions $f:\bb{R}^d\times\Omega\longrightarrow \bb{R}$ with 
%$f(\cdot,\omega)\in C^k(\bb{R}^d)$ for each $\omega\in\Omega$ and such that for all 
%$x,y\in\bb{R}^d,\ \omega\in\Omega$ and $J\in\{0,1,\ldots,k\}^d$ with 
%$|J|\leq k,$\begin{equation}\label{20}
%|\partial^Jf(x,\omega)|\leq m,\quad
% |\partial^Jf(x,\omega)-
%\partial^Jf(y,\omega)|
%\leq m|x-y|.
%\end{equation}
For a function \begin{equation}\label{301}
\varphi\in \mathrm{Lip}^2_1(\bb{R}^d\times\Omega)
\end{equation}
 which satisfies the finite range dependence condition with range $R/2,$ i.e. for all 
Borel subsets $A,B\subseteq\bb{R}^d,$\begin{equation}\label{302}
\cal{H}_A^{\varphi}\mbox{ and }\cal{H}_B^{\varphi}\mbox{ are $\bb{P}$-independent whenever }d(A,B)>R/2,
\end{equation}where $\cal{H}_A^{\varphi},\cal{H}_B^{\varphi}$ are defined in (\ref{300}), we define for $0\leq \eps\leq 1,\ x\in\bb{R}^d$
and $\omega\in\Omega,$\begin{equation}\label{22}
\phi_{\eps}(x,\omega)\df\frac{1+\frac{\eps}{d}\varphi(x,\omega)}
{1+\frac{\eps}{d}E^{\bb{P}}[\varphi(0,\omega)]}.
\end{equation}Let us collect some useful properties of the family $\phi_{\eps},\ 0\leq \eps\leq 1,$
in the following
{\lemma{\label{l1}
For every Borel subset $A\subseteq\bb{R}^d$ we have that for all $0\leq \eps\leq 1,$
\begin{equation}\label{101}
\cal{H}^{\phi_{\eps}}_{A} \subseteq \cal{H}^{\varphi}_{A}.
\end{equation}
Moreover, the following 
properties hold true for all $0\leq \eps\leq 1,\ x,y\in\bb{R}^d,\ \omega\in\Omega$ and 
$i=1,\ldots,d:$ (note that $d\geq 2$) 
\begin{eqnarray}
\label{103}&&\phi_{\eps}\in\mathrm{Lip}^2_{3}(\bb{R}^d\times\Omega);\\
\label{71}&&\frac{1}{3}\leq \frac{d-\eps}{d+\eps}\leq \phi_{\eps}(x,\omega)\leq \frac{d+\eps}{d-\eps}
\leq 3;\\
\label{72}&&|\partial_i^j\phi_{\eps}(x,\omega)|\leq\frac{\eps}{d-\eps}\leq 1\quad\mbox{ for }
j=1,2;\\
\label{73}&&|\partial_i^j\phi_{\eps}(x,\omega)-\partial_i^j\phi_{\eps}(y,\omega)|\leq\frac{\eps}{d-\eps}|x-y|
\quad\mbox{ for } j=0,1,2;\\
\label{74}&&d\bb{Q}_{\eps}\df \phi_{\eps}(0,\omega)d\mathbb{P}\mbox{ defines a probability 
measure equivalent to }\bb{P}.
\end{eqnarray}}}
%{\lemma{\label{l1}
%For a Borel subset $A\subseteq\bb{R}^d$ we have that for all $0\leq \eps\leq 1,$
%\begin{equation}\label{101}
%\cal{H}^{\phi_{\eps}}_{A} \subseteq \cal{H}^{\varphi}_{A}.
%\end{equation}
%Moreover, the following 
%properties hold true for all $0\leq \eps\leq 1,\ x,y\in\bb{R}^d,\ \omega\in\Omega$ and 
%$J\in\{0,1,2\}^d:$ (note that $d\geq 2$) 
%\begin{eqnarray}
%\label{103}&&\phi_{\eps}\in\mathrm{Lip}^2_{\Phi}(\bb{R}^d\times\Omega)\mbox{ with }\Phi=3;\\
%\label{71}&&\frac{1}{3}\leq \frac{d-\eps}{d+\eps}\leq \phi_{\eps}(x,\omega)\leq \frac{d+\eps}{d-\eps}
%\leq 3;\\
%\label{72}&&|\partial^J\phi_{\eps}(x,\omega)|\leq\frac{\eps}{d-\eps}\leq
% 1\quad\mbox{ for }
%1\leq |J|\leq 2;\\
%\label{73}&&|\partial^J\phi_{\eps}(x,\omega)-\partial^J
%\phi_{\eps}(y,\omega)|\leq\frac{\eps}{d-\eps}|x-y|
%\quad\mbox{ for } 0\leq |J|\leq 2;\\
%\label{74}&&\bb{Q}_{\eps}[\,\star\,]\df E^{\bb{P}}[\,\star\, ,\phi_{\eps}(0,\omega)]\mbox{ defines a probability 
%measure equivalent to }\bb{P}.
%\end{eqnarray}}}

{\textbf{Proof:}} The measurability property (\ref{101}) is obvious. Due to assumption (\ref{301}) we find that  
$\phi_{\eps}(\cdot,\omega)\in C^2(\bb{R}^d)$ for all $0\leq \eps\leq 1,\ 
\omega\in\Omega.$ The statements (\ref{71})-(\ref{73}) follow by direct inspection of the formula (\ref{22}) and 
property (\ref{74}) is a direct consequence of (\ref{71}) and the fact that $E^{\bb{P}}[\phi_{\eps}(0,\omega)]=1$
for all $0\leq \eps\leq 1.$ Finally observe that for all  $0\leq \eps\leq 1,\ 
\phi_{\eps}\in \mathrm{Lip}^2_{3}(\bb{R}^d\times\Omega).$ 
\begin{flushright}$\Box$\end{flushright}

Furthermore, we consider a measurable function $h(x,\omega)=(h_{i,j}(x,\omega))_{i,j=1,\ldots,d},\ x\in\bb{R}^d,
\ \omega\in\Omega,$ with values in the space of skew-symmetric $d\times d$-matrices, i.e.
\begin{equation}\label{303}
h_{i,j}(x,\omega)=-h_{j,i}(x,\omega)\qquad\mbox{for all }x\in\bb{R}^d,\omega\in\Omega\mbox{ and }
i,j=1,\ldots,d,
\end{equation}and we assume that for all $i,j=1,\ldots,d,$ \begin{equation}\label{75}
h_{i,j}\in {\mathrm{Lip}}^2_1(\bb{R}^d\times\Omega)\qquad\mbox{and}\qquad 
\cal{H}^{h_{i,j}}_{A}\subseteq\cal{H}^{\varphi}_{A^{R/8}}
\end{equation}for all Borel subsets $A$ of $\bb{R}^d.$
We then define the function $c=(c_1,\ldots,c_d)^T:\bb{R}^d\times\Omega\longrightarrow\bb{R}^d$
as\begin{equation}\label{24}
c_i(x,\omega)\df\frac{1}{8d^2}\sum_{j=1}^{d}\partial_j h_{i,j}(x,\omega),\quad\mbox{ for all }
x\in \bb{R}^d,\omega\in\Omega\mbox{ and }i=1,\ldots,d.
\end{equation}By direct inspection of the definition (\ref{24}) and using the properties mentioned in (\ref{75}) one easily finds 
{\lemma{\label{l2}
For all $x,y\in\bb{R}^d$ and 
$\omega\in\Omega$ we have that $c_i(\cdot,\omega)\in C^1(\bb{R}^d),\ i=1,\ldots,d,$ and 
\begin{equation}\label{25}
|c(x,\omega)|\leq \frac{1}{8},\quad |c(x,\omega)-c(y,\omega)|\leq \frac{1}{8}|x-y|,\quad 
\nabla\cdot c(x,\omega)=0.
\end{equation}Moreover, the following measurability property holds: for every Borel subset $A$ of $\bb{R}^d,$
\begin{equation}\label{77}
\cal{H}^c_{A} \subseteq\cal{H}^{\varphi}_{A^{R/4}}.
\end{equation}}}
%{\lemma{\label{l2}
%For $i=1,\ldots,d,\ c_i(\cdot,\omega)\in \mathrm{Lip}^1_{1/8d}
%(\bb{R}^d\times\Omega)$ and for all $x,y\in\bb{R}^d$ and 
%$\omega\in\Omega$ we have that
%\begin{equation}\label{25}
%|c(x,\omega)|\leq \frac{1}{8},\quad |c(x,\omega)-c(y,\omega)|\leq \frac{1}{8}|x-y|,\quad 
%\nabla\cdot c(x,\omega)=0.
%\end{equation}Moreover, the following measurability property holds: for a Borel subset $A$ of $\bb{R}^d,$
%\begin{equation}\label{77}
%\cal{H}^c_{A} \subseteq\cal{H}^{\varphi}_{A^{R/4}}.
%\end{equation}}}
%{\rem{\label{26}In stead of assuming that all first derivatives with respect to $x\in\bb{R}^d$of the 
%entries $h_{i,j}(x,\omega)$ are bounded and Lipschitz continuous with constant 1, it suffices to
% impose the following slightly weaker conditions to $h_{i,j}$: for all $x,y\in\bb{R}^d,\ \omega\in
% \Omega$ and $i,j=1,\ldots, d,$\begin{equation}
% |\partial_j h_{i,j}(x,\omega)|\leq 1,\qquad |\partial_j h_{i,j}(x,\omega)-\partial_j h_{i,j}
% (y,\omega)|\leq |x-y|.
% \end{equation}}}

In addition to the above restrictions on the choices of $c$ and $\varphi$ we assume the following 
non-degeneracy condition to be satisfied:\begin{equation}\label{41}
E^{\bb{P}}\left[ c(0,\omega)\varphi(0,\omega)\right]\ne 0 .
\end{equation}
The reason for this constraint will become clear when we construct the first example (see (\ref{100})). 
Note that since for all $i,j=1,\ldots, d,\ h_{i,j}\in\mathrm{Lip}^2_1(\bb{R}^d\times\Omega),$ both 
$h_{i,j}(x,\omega)$ and 
$\partial_j h_{i,j}(x,\omega)$ are bounded in absolute value by 1 for all $x\in\bb{R}^d,\omega\in\Omega$ and hence 
$E^{\bb{P}}[c_i(0,\omega)]\overset{{\scriptscriptstyle{(\ref{24})}}}{=}
(8d^2)^{-1}\sum_{j=1}^d E^{\bb{P}}[\partial_j h_{i,j}(0,\omega)]=(8d^2)^{-1}
\sum_{j=1}^d \partial_j E^{\bb{P}}[h_{i,j}(0,\omega)].$ Due to stationarity of the environment the latter expression 
equals zero. So, we find that\begin{equation}\label{110}
E^{\bb{P}}\left[c(0,\omega)\right]=0.
\end{equation}This together with (\ref{41}) implies that\begin{equation}\label{78}
\varphi(0,\omega)\not\equiv \mathrm{const.}\qquad \bb{P}\mbox{-a.s..}
\end{equation}
For a possible example of functions $\varphi$ and $h_{i,j}, i,j=1,\ldots d$ with the properties mentioned in 
(\ref{301}), (\ref{302}), and  (\ref{303}), (\ref{75}) respectively, such that (\ref{41}) holds, we refer to Section \ref{e2}.
Now we are ready to introduce the family of drifts announced in the introduction of this section.
In the notation of (\ref{22}) and 
(\ref{24}), for $0\leq\eps\leq 1$ and  $\lambda\in [-1,1],$ we define the drift\begin{equation}
\label{14}
b_{\varepsilon,\lambda}(x,\omega)\df\frac{\nabla\phi_{\varepsilon}(x,\omega)}
{2\phi_{\varepsilon}(x,\omega)}+
\eps\left( \frac{ c(x,\omega)+\lambda E^{\bb{P}}\left[ \frac{c(0,\omega)}{\phi_{\eps}(0,\omega)}\right]}
{\phi_{\varepsilon}(x,\omega)}\right),\qquad x\in\bb{R}^d,\omega\in\Omega.
\end{equation}
{\lemma{\label{l3}
There is an $0< \eps_0\leq 1$ such that for all $0\leq\eps\leq\eps_0$ and $\lambda\in [-1,1],$ the drift $b_{\eps,\lambda}$ defined in (\ref{14}) 
is finite range dependent with range $R$ and for all $x,y\in\bb{R}^d, \omega\in\Omega,$ the 
following holds:\begin{equation}\label{27}
|b_{\eps,\lambda}(x,\omega)|\leq \eps,\quad |b_{\eps,\lambda}(x,\omega)-b_{\eps,\lambda}(y,\omega)|
\leq |x-y|.
\end{equation}}}

The proof of Lemma \ref{l3} is given in Section \ref{e1}.\\

For $0\leq \eps\leq \eps_0, \lambda\in [-1,1], x\in\bb{R}^d$ and $\omega\in\Omega$ we denote with $P_{x,\omega}^{\eps,\lambda}$
the unique solution of the martingale problem starting in $x$ at time 0 and attached to \begin{equation}\label{13}
\cal{L}^{\omega,\varepsilon,\lambda}=\frac{1}{2}\Delta+ b_{\varepsilon,\lambda}(\cdot,\omega)\cdot\nabla,
\end{equation}see Corollary 5.4.29 in \cite{ks} for the 
existence of a unique solution of the above martingale problem. According to the definition (\ref{6}) we denote
by $P_{0}^{\eps,\lambda}$ the {\textit{annealed}} law $\bb{P}\times P^{\eps,\lambda}_{0,\omega}.$ We know that the process of the 
environment viewed from the particle associated to the diffusion $P_{0,\omega}^{\eps,\lambda},$ which is defined in analogy to
(\ref{11}), has $\bb{Q}_{\varepsilon}$ (see (\ref{74})) as an invariant ergodic measure
 equivalent to $\bb{P}$ and
\begin{equation}\label{17}P_0^{\eps,\lambda}\mbox{-a.s.,}\quad 
\frac{X_t}{t}\overset{t\to\infty}{\longrightarrow}
v_{\eps,\lambda}\df E^{\bb{Q}_{\varepsilon}}\left[b_{\eps,\lambda}(0,\omega)\right],
\end{equation}see Theorem \ref{t01} in the Appendix. For the expected local drift under the static measure $\bb{P}$ 
we write\begin{equation}
\label{18}
d_{\eps,\lambda}\df E^{\bb{P}}\left[b_{\eps,\lambda}(0,\omega)\right].
\end{equation}

\section{The Examples}\label{c}

In this section we are going to provide examples of diffusions in random environment with arbitrarily small, non-vanishing 
expected local drift but without ballistic behavior. A slight modification of this family of diffusions 
gives us then examples of ballistic 
behavior, but where the expected local drift under the static measure 
vanishes or even has an opposite direction to the limiting velocity.\\

In order to shorten notation let us introduce the 
continuous function\begin{eqnarray}\label{31}
\label{30} g:\left[0,1\right] \longrightarrow  \left[1,3\right];\quad
 \eps \mapsto  g(\eps)\df E^{\bb{P}}\left[ \phi_{\eps}
(0,\omega)^{-1}\right]\overset{\mathrm{Jensen}}{\geq} E^{\bb{P}}\left[ \phi_{\eps}
(0,\omega)\right]^{-1}\overset{(\ref{74})}{=}1,
\end{eqnarray}where $g\leq 3$ follows from the uniform lower bound on $\phi_{\eps}$ given in (\ref{71}).\\

%Since $g(\cdot)$ is continuous, the supremum on a compact interval is attained and therefore 
%the following quantity is well-defined:
%\begin{equation}\label{90}
%\delta(\eps)\df\max_{\eta\in [0,\eps]}g(\eta),\qquad \mbox{for }0\leq\eps\leq 1/10.
%\end{equation} 

Now we are ready to provide our first examples. Recall the definition of $\eps_0$ in Lemma \ref{l3}.
{\thm{\label{p1}
($\lambda=0$).  
There exists an $0<\eps_1\leq \eps_0$ such that for all $0<\eps<\eps_1$ we have that
 $v_{\eps,0}=0$ but $d_{\eps,0}\ne 0.$
}}\\

{\textbf{Proof:}} Let us begin with the calculation of the limiting velocity. Let $\lambda=0$ and 
$0\leq\eps\leq \eps_0.$
Then\begin{eqnarray}
\label{112} v_{\eps,0}&
\overset{(\ref{17})}{=}&E^{\bb{Q}_{\eps}}\left[b_{\eps,0}(0,\omega)\right]
\overset{(\ref{74}),(\ref{14})}{=}E^{\bb{P}}\left[\frac{1}{2}\nabla\phi_{\eps}(0,\omega)\right]
+\eps E^{\bb{P}}\left[c(0,\omega)\right]=0,
\end{eqnarray}where in the last equality we used (\ref{110}) and that $E^{\bb{P}}\left[\nabla\phi_{\eps}(0,\omega)\right]=
\nabla E^{\bb{P}}\left[\phi_{\eps}(0,\omega)\right]=0,$ which can be shown by arguments similar to those in the 
derivation of (\ref{110}). For the local drift under the static measure we have:
\begin{equation*}
d_{\eps,0}\overset{(\ref{18})}{=} E^{\bb{P}}\left[b_{\eps,0}(0,\omega)\right]
\overset{(\ref{14})}{=}\frac{1}{2}E^{\bb{P}}\left[\nabla \log\left(\phi_{\eps}(0,\omega)\right)\right]
+\eps E^{\bb{P}}\left[\frac{c(0,\omega)}{\phi_{\eps}(0,\omega)}\right].
\end{equation*}Again by arguments similar to those leading to (\ref{110}) we find that the first term on the right-hand side of the second equality above vanishes and hence
\begin{equation}\label{111}
d_{\eps,0}=\eps E^{\bb{P}}\left[\frac{c(0,\omega)}{\phi_{\eps}(0,\omega)}\right]
\overset{(\ref{22})}{=} \eps \left(1+\frac{\eps}{d} E^{\bb{P}}\left[\varphi(0,\omega)\right]\right)
E^{\bb{P}}\left[\frac{c(0,\omega)}{1+\frac{\eps}{d}\varphi(0,\omega)}\right].
\end{equation}

%&=&\frac{1}{2}\nabla E^{\bb{P}}\left[\log\left(\phi_{\eps}(0,\omega)\right)\right]+
%\eps E^{\bb{P}}\left[\frac{c(0,\omega)}{\phi_{\eps}(0,\omega)}\right]\\
%&\overset{stat.}{=}&\eps E^{\bb{P}}\left[\frac{c(0,\omega)}{\phi_{\eps}(0,\omega)}\right]=
%\eps \left(1+\frac{\eps}{d} E^{\bb{P}}\left[\varphi(0,\omega)\right]\right)
%E^{\bb{P}}\left[\frac{c(0,\omega)}{1+\frac{\eps}{d}\varphi(0,\omega)}\right].\\
Using the identity (\ref{111}) one can check by straightforward computations that 
%\begin{equation*}
%d_{0,0}=0,\quad \frac{d}{d\eps}\Big|_{\eps=0}d_{\eps,0}=E^{\bb{P}}[c(0,\omega)]\overset{(\ref{110})}{=}0
%\end{equation*}and
\begin{eqnarray}
\nonumber d_{0,0}=0,&& \frac{d}{d\eps}\Big|_{\eps=0}d_{\eps,0}=E^{\bb{P}}[c(0,\omega)]\overset{(\ref{110})}{=}0
\\
\nonumber \mbox{and}\qquad\qquad\qquad\qquad\qquad\qquad&& \\
\nonumber 
\frac{d^2}{d\eps^2}\Big|_{\eps=0}d_{\eps,0}&=&\frac{2}{d}E^{\bb{P}}\left[\varphi(0,\omega)\right]
E^{\bb{P}}\left[c(0,\omega)\right]-\frac{2}{d}E^{\bb{P}}\left[\varphi(0,\omega)c(0,\omega)\right]\\
\label{100}&\overset{(\ref{110})}{=}&-\frac{2}{d}E^{\bb{P}}\left[\varphi(0,\omega)c(0,\omega)\right]
\overset{(\ref{41})}{\ne}0.
\end{eqnarray}Since $d_{\eps,0}$ is continuous in $\eps\in [0,1],$ it then follows from the above computations that there exists an $0<\eps_1\leq \eps_0$ such that for all $0<\eps<\eps_1,
\ d_{\eps,0}\ne 0$ but at the same time we have that $v_{\eps,0}=0,$ as shown in (\ref{112}).
\begin{flushright}$\Box$\end{flushright}

Recall the definition of $\eps_1$ in Theorem \ref{p1}.

 {\thm{\label{p2}
  For all 
 $0<\eps<\eps_1$ and $\lambda\in (-1/3,0)$ 
 there exists a constant $\gamma>0$ depending on $\eps$ and $\lambda$ such that 
 $d_{\eps,\lambda}=-\gamma v_{\eps,\lambda}\ne 0.$
 }}\\
 
 {\textbf{Proof:}} Recall the definition of $g$ given in (\ref{31}). As a consequence of Theorem \ref{p1} we know that for all 
 $0<\eps<\eps_1$ and $\lambda\in [-1,1]\backslash\{0\},$ the 
 limiting velocity is not equal to zero. Indeed,\begin{equation}\label{33}
 v_{\eps,\lambda}\overset{(\ref{17})}{=}E^{\bb{Q}_{\eps}}\left[ b_{\eps,\lambda}(0,\omega)\right]
 \overset{(\ref{14}),(\ref{111})}{=}
 v_{\eps,0}+\lambda d_{\eps,0}\overset{(\ref{112})}{=}\lambda d_{\eps,0}\ne 0. 
 \end{equation}Furthermore, for all $\lambda\in (-1/3,0)$ and 
 $0<\eps\leq \eps_0,$ see Lemma \ref{l3} for the definition of $\eps_0,$ the expected local drift under the 
 static measure $\bb{P}$ equals\begin{equation}\label{34}
 d_{\eps,\lambda}\overset{(\ref{18}),(\ref{14}),(\ref{111})}{=}d_{\eps,0}\left( 1+\lambda g(\eps)\right)
 \overset{(\ref{33})}{=}
 \lambda^{-1}v_{\eps,\lambda} \left( 1+\lambda g(\eps)\right)=
 -\gamma v_{\eps,\lambda}
 \end{equation}with $\gamma\df-\lambda^{-1}( 1+\lambda g(\eps))>0,$ where we used that $\lambda\in (-1/3,0)$ and  $g\leq 3.$  \begin{flushright}$\Box$\end{flushright}
{\rem{$\mathrm{W}$}}e know that in the one-dimensional discrete setting, where we consider random walks in an i.i.d. 
random environment, 
the limiting velocity cannot have an opposite direction to the expected local drift under the static measure, see 
Remark on page 211 in \cite{sznit3}. In the case of one-dimensional diffusions in random environment 
as described in the 
Introduction, i.e. diffusions generated by operators of the form (\ref{5}), the same holds true. Indeed, by Proposition 
2.7 and formula (2.77) in \cite{goer} we know that\begin{equation*}
E^{\bb{P}}\left[b(0,\omega)\right]>0\Longleftrightarrow P_0\mbox{-a.s., }\lim_{t\to\infty}X_t=+\infty,
\end{equation*} which implies that \begin{equation*}
P_0\mbox{-a.s., }\liminf_{t\to\infty}\frac{X_t}{t}\geq 0.
\end{equation*}\begin{flushright}$\Box$\end{flushright}%\pagebreak
In the next theorem we will see that by examining more carefully the formula for the expected drift under the
static measure $d_{\eps,\lambda}{=}d_{\eps,0}\left( 1+\lambda g(\eps)\right),$  which was derived in (\ref{34}), we can find examples of diffusions with vanishing expected local drift but still with 
ballistic behavior. Recall the definition of $\eps_1$ in Theorem \ref{p1}.

{\thm{\label{p3}
 For all $0<\eps<\eps_1$ we find a $\lambda\in [-1,-1/3]$ 
depending on $\eps$ such that $v_{\eps,\lambda}\ne 0$ and $d_{\eps,\lambda}=0.$ 
}}\\

{\textbf{Proof:}} Recall the definition of $g$ given in (\ref{31}). Let $0<\eps<\eps_1$ and choose $\lambda\df -g(\eps)^{-1}\in [-1,-1/3].$ Then 
$v_{\eps,\lambda}\ne 0,$ which is shown in (\ref{33}), and $d_{\eps,\lambda}=d_{\eps,0}(1+\lambda g(\eps))=0.$\begin{flushright}$\Box$\end{flushright}

%{\prop{\label{p3}
%Recall the definition of $\eps_1$ in Proposition \ref{p1} and $\delta(\cdot)$ in (\ref{90}). 
%For all $\eps_2\in (0,\eps_1)$ and all $\lambda\in (-1,-1/\delta(\eps_2)]\subseteq (-1,-19/21]$ there exists an 
%$\eps\in (0,\eps_2]$ depending on $\lambda$ such that $v_{\eps,\lambda}\ne 0$
% and $d_{\eps,\lambda}=0.$ 
%}}\\

%{\textbf{Proof:}} Let $\eps_2\in (0,\eps_1)$ and pick a 
%$\lambda\in (-1,-1/\delta(\eps_2)].$ Then for all $\eps\in (0,\eps_2]\subseteq (0,\eps_1),$
%$v_{\eps,\lambda}{\ne} 0\mbox{ and } d_{\eps,\lambda}{=}
%d_{\eps,0}\left( 1+\lambda g(\eps)\right),$
%see (\ref{33}) and (\ref{34}). From the definition of $\delta,$ see (\ref{90}), the 
%continuity of $g$ and the fact that $g(0)=1,$ we know that there exist
%${\eps^{*}}\in (0,\eps_2]$ 
%which depend on $\lambda$ such that $g(\eps^{*})=|\lambda|^{-1}\in (1,\delta(\eps_2)]$ and hence 
%$1+\lambda g(\eps^{*})=1-|\lambda||\lambda|^{-1}=0$ which implies that $d_{\eps^*,\lambda}=0.$ 
%This finishes the proof.
%\begin{flushright}$\Box$\end{flushright}

\section{Auxiliary results}\label{e}

In this section we give the proof of Lemma \ref{l3} and provide a concrete example of functions 
$\varphi$ and $h_{i,j},i,j=1,\ldots d,$
with the properties given in (\ref{301}), (\ref{302}), and (\ref{303}), (\ref{75}) respectively, such that (\ref{41}) holds. These 
functions are involved in the construction of the family of drifts defined in (\ref{14}).

\subsection{Proof of Lemma \ref{l3}}\label{e1}
It follows from the definition (\ref{14}) and the measurability properties (\ref{101})
and (\ref{77}) that for all $0\leq\eps\leq 1$ and $\lambda\in[-1,1],\ 
\cal{H}^{b_{\eps,\lambda}}_{A} \subseteq\cal{H}^{\varphi}_{A^{R/4}},$
and hence the drift $b_{\eps,\lambda}$ satisfies the finite range dependence condition with range $R$
due to the finite range dependence of $\varphi$ with range $R/2$ mentioned in (\ref{302}).
With the help of Lemma \ref{l1} and \ref{l2} 
%and the fact that\begin{equation}\label{81}
%\frac{d}{d-\eps}\leq \frac{20}{19}\qquad\mbox{ for all }0\leq\eps\leq 1/10,
%\end{equation}
we find that for all 
$0\leq\eps\leq 1,\lambda\in[-1,1], x\in\bb{R}^d$ and $\omega\in\Omega,$\begin{eqnarray*}
|b_{\eps,\lambda}(x,\omega)|&\overset{(\ref{71})}{\leq} &\frac{d+\eps}{d-\eps}\left\lbrace 
\frac{1}{2}|\nabla\phi_{\eps}{(x,\omega)}|+\eps |c(x,\omega)|+\frac{d+\eps}{d-\eps}\eps|\lambda|
E^{\bb{P}}[|c(0,\omega)|]\right\rbrace \\
&\overset{(\ref{72}),(\ref{25})}{\leq}&\eps\frac{d+\eps}{d-\eps}\left\lbrace 
\frac{1}{2}\frac{d}{d-\eps}+\frac{1}{8}+\frac{1}{8}\frac{d+\eps}{d-\eps}\right\rbrace 
\leq \eps,
\end{eqnarray*}where the last inequality holds for all $0\leq\eps\leq\delta_1$ with $\delta_1>0$ small enough.
It remains to prove the Lipschitz property in (\ref{27}).
%$0\leq \eps\\leq 1/10,\lambda\in[-1,1], x\in\bb{R}^d$ and $\omega\in\Omega,$\begin{eqnarray*}
%|b_{\eps,\lambda}(x,\omega)|&\overset{(\ref{71})}{\leq} &\frac{21}{19}\left\lbrace 
%\frac{1}{2}|\nabla\phi_{\eps}{(x,\omega)}|+\eps |c(x,\omega)|+\frac{21}{19}\eps|\lambda|
%E^{\bb{P}}[|c(0,\omega)|]\right\rbrace \\
%&\overset{(\ref{72}),(\ref{25})}{\leq}&\frac{21}{19}\left\lbrace 
%\frac{1}{2}d\frac{\eps}{d-\eps}+\eps\frac{1}{8}+\frac{21}{19\cdot8}\eps|\lambda|\right\rbrace 
%\leq \frac{315}{361}\eps,
%\end{eqnarray*}where in the last inequality we used that $|\lambda|\leq 1$ and 
%This finishes the proof of the first bound in (\ref{27}).
For the remainder of the proof let us define for all 
$0\leq\eps\leq 1$ and $\lambda\in [-1,1],$\begin{equation}\label{79}
\alpha_{\eps,\lambda}\df\lambda E^{\bb{P}}\left[ \frac{c(0,\omega)}{\phi_{\eps}(0,\omega)}\right]
\overset{(\ref{71}),(\ref{25})}{\leq}\vartheta_0,
\end{equation}for some constant $\vartheta_0>0.$ 
For $\eps, \lambda$ as above, $x,y$ in $\bb{R}^d$ and $\omega$ in $\Omega,$ we obtain\begin{eqnarray}
\nonumber &&|b_{\eps,\lambda}(x,\omega)-b_{\eps,\lambda}(y,\omega)|\\
\nonumber&\leq&\frac{1}{2\phi_{\eps}(x,\omega)\phi_{\eps}(y,\omega)}\bigg\{
\left|\nabla\phi_{\eps}(x,\omega)\phi_{\eps}(y,\omega)-\nabla\phi_{\eps}(y,\omega)\phi_{\eps}(x,\omega)\right|\\[5pt]
\nonumber &&+2\eps\left|c(x,\omega)\phi_{\eps}(y,\omega)-c(y,\omega)\phi_{\eps}(x,\omega)\right|
+2\eps\left|\alpha_{\eps,\lambda}\right|\left|\phi_{\eps}(y,\omega)-\phi_{\eps}(x,\omega)\right|\bigg\}\\
\nonumber &\df&\frac{f_{\eps}(x,y,\omega)+2\eps g_{\eps}
(x,y,\omega)+2\eps|\alpha_{\eps,\lambda}|h_{\eps}(x,y,\omega)}{2\phi_{\eps}(x,\omega)\phi_{\eps}(y,\omega)}\\
\label{102}&\overset{(\ref{71})}{\leq}& \frac{9}{2}\left(f_{\eps}(x,y,\omega)+2\eps g_{\eps}
(x,y,\omega)+2\eps|\alpha_{\eps,\lambda}|h_{\eps}(x,y,\omega)\right).
\end{eqnarray}
We also find that $f_{\eps}(x,y,\omega)$ is smaller or equal to
\begin{eqnarray}
%&&\nonumber\hspace*{0.6cm}\left|\nabla\phi_{\eps}(x,\omega)\phi_{\eps}(y,\omega)-\nabla\phi_{\eps}(x,\omega)
%\phi_{\eps}(x,\omega)\right|+
%\left|\nabla\phi_{\eps}(x,\omega)\phi_{\eps}(x,\omega)-\nabla\phi_{\eps}(y,\omega)\phi_{\eps}(x,\omega)\right|\\[5pt]
&&\nonumber\left|\nabla\phi_{\eps}(x,\omega)\right|
\left|\phi_{\eps}(y,\omega)-\phi_{\eps}(x,\omega)\right|+\left|\phi_{\eps}(x,\omega)\right|
\left|\nabla\phi_{\eps}(x,\omega)-\nabla\phi_{\eps}(y,\omega)\right|\\
&&\label{82}\overset{(\ref{71})-(\ref{73})}{\leq}d\left(\frac{\eps}{d-\eps}\right)^2|x-y|+\frac{d+\eps}{d-\eps} \cdot\frac{d\eps}{d-\eps}|x-y|\leq
\eps\vartheta_1|x-y|,
\end{eqnarray}for a positive constant $\vartheta_1,$ which depends on the upper bounds given in
(\ref{71}), (\ref{72}). By an analogous
computation using (\ref{71})-(\ref{73}) and (\ref{25}) one can show that for some constant $\vartheta_2>0,$ \begin{equation}\label{83}
g_{\eps}(x,y,\omega)\leq \vartheta_2 |x-y|.
\end{equation}Finally, from (\ref{72}) and (\ref{73}) follows that there is a 
constant $\vartheta_3>0$ such that\begin{equation}\label{84}
h_{\eps}(x,y,\omega)\leq \vartheta_3|x-y|.
\end{equation}Collecting (\ref{79}) - (\ref{84}), we find the Lipschitz property (\ref{27}) for all 
$0\leq\eps\leq\delta_2$ with $\delta_2>0$ small enough. Our claim then follows with $\eps_0$ equal to the 
minimum of $\delta_1$ and $\delta_2.$ This finishes the proof of Lemma \ref{l3}.

\subsection{A concrete example}\label{e2}
A possible example of functions $\varphi$ and $h_{i,j}, i,j=1,\ldots d$ with the properties mentioned in 
(\ref{301}), (\ref{302}), and  (\ref{303}), (\ref{75}) respectively, such that (\ref{41}) holds,
can be constructed as follows. As a random 
environment $(\Omega,\cal{A},\bb{P})$ we consider a canonical Poisson point process on 
$\bb{R}^d$ with constant intensity and for all 
$x\in\bb{R}^d,\omega\in\Omega$ and a Borel subset $A\subseteq\bb{R}^d$ we define the group 
of transformation on $\Omega$ as $\tau_x(\omega)(A)\df\omega(x+A),$ where 
$x+A\df\{x+a\,|\,a\in A\}.$ Assume for 
the moment that we have a function $\varphi\in \mathrm{Lip}^3_1(\bb{R}^d\times\Omega)$
with $\partial_i\varphi\in \mathrm{Lip}^2_1(\bb{R}^d\times\Omega)$ satisfying the
finite range dependence condition with range $R/2$ and (\ref{78}). A possible skew-symmetric matrix 
$(h_{i,j}(x,\omega))_{i,j=1,\ldots,d},\ x\in\bb{R}^d,\ \omega\in\Omega$ 
is defined as follows:\begin{eqnarray*}
&& h_{i,i}(x,\omega)=0,\quad \mbox{for }i=1,\ldots,d;\\
&& h_{i,1}(x,\omega)=\partial_i \varphi(x,\omega)\quad\mbox{and}\quad
h_{1,i}(x,\omega)=-\partial_i \varphi(x,\omega),\quad \mbox{for }i=2,\ldots,d;\\
&&(h_{i,j}(x,\omega)),\ i,j=2,\dots,d,\mbox{ is a }(d-1)\times (d-1)\mbox{-dimensional skew-symmetric 
ma-}\\
&&\mbox{trix with entries which fulfill (\ref{75}).}
\end{eqnarray*}
One then easily sees that the entries $h=(h_{i,j}),\ i,j=1,\ldots,d,$ satisfy (\ref{75}) and hence the function $c$ defined via the formula (\ref{24}) 
fulfills the properties in Lemma \ref{l2}. The following lines show that condition (\ref{41}) holds as well.
Indeed assume the contrary, then: \begin{eqnarray*}
0&=& E^{\bb{P}}\left[c_1(0,\omega)\varphi(0,\omega) \right]=
-\frac{1}{8d^2}\sum_{j=2}^d E^{\bb{P}}\left[\left( \partial^2_j\varphi(0,\omega)\right) \varphi(0,\omega) \right]\\
&=&-\frac{1}{8d^2}\sum_{j=2}^d E^{\bb{P}}\left[\partial_j\left( \partial_j\varphi(0,\omega) \varphi(0,\omega) \right)\right]
+\frac{1}{8d^2}\sum_{j=2}^d E^{\bb{P}}\left[\left( \partial_j\varphi(0,\omega)\right)^2\right]\\
&{=}&
-\frac{1}{8d^2}\sum_{j=2}^d \partial_j E^{\bb{P}}\left[\partial_j\varphi(0,\omega) \varphi(0,\omega) \right]+
\frac{1}{8d^2}\sum_{j=2}^d E^{\bb{P}}\left[\left( \partial_j\varphi(0,\omega) \right)^2\right]\\
&=&\frac{1}{8d^2}\sum_{j=2}^d E^{\bb{P}}\left[\left( \partial_j\varphi(0,\omega) \right)^2\right],
\end{eqnarray*}
where we used stationarity in the last 
equality. This implies that for $\bb{P}$-a.e. $\omega\in\Omega,$\begin{equation}
\partial_j\varphi(0,\omega)=0,\mbox{ for all }j=2,\ldots,d,
\end{equation}and hence by stationarity and continuity of $\partial_j\varphi(\cdot,\omega)$ for each $\omega\in\Omega,$ 
we find that $\bb{P}$-a.s.,\begin{equation}
\partial_j\varphi(x,\omega)=0,\mbox{ for all }x\in\bb{R}^d\mbox{ and }j=2,\ldots,d.
\end{equation}It follows that for $\bb{P}$-a.e. $\omega\in\Omega,\ \varphi(x,\omega)$ is in fact 
a function of $x_1$ and $\omega$ only. In the notation $\bar{R}_2=(0,R,0,\ldots,0)\in\bb{R}^d$ we thus have that 
$\bb{P}\mbox{-a.s., }\varphi(0,\omega)=\varphi(\bar{R}_2,\omega)$ and since $\varphi$ satisfies the finite 
range condition with range $R/2$ we find that for all integers $n\geq 0,$
\begin{equation*}
E^{\bb{P}}\left[\varphi(0,\omega)^n\right]{=}
E^{\bb{P}}\left[\prod_{k=0}^{n-1}\varphi(k\bar{R}_2,\omega)\right]\overset{indep.}{=}
\prod_{k=0}^{n-1}E^{\bb{P}}\left[\varphi(k\bar{R}_2,\omega)\right]{=}
E^{\bb{P}}\left[\varphi(0,\omega)\right]^n
\end{equation*}which is a contradiction to (\ref{78}) and (\ref{41}) must hold true.\\
We now come to the construction of a possible 
$\varphi$ satisfying the above required conditions. Pick a 
measurable real-valued non-negative function $\zeta(x),\ x\in\bb{R}^d,$ which is supported in a ball of radius $R/8$ 
and strictly positive on a set of positive Lebesgue measure. 
 Then convolve the Poisson point process with
$\zeta$ and truncate, let us say at 1, the new function, i.e. for $x\in\mathbb{R}^{d},\omega\in\Omega$ we define\begin{equation} 
\hat{\varphi}(x,\omega)\df\left(\int_{\mathbb{R}^d}\zeta(x-y)\omega(dy)\right)\wedge 1.
                                           \end{equation}
After smoothing with a non-negative  
mollifier $\rho$ belonging to $C^3(\bb{R}^d)$ which is also supported in a ball of radius $R/8$ we obtain 
a function\begin{equation} 
\tilde{\varphi}(x,\omega)\df\int_{\mathbb{R}^d}\hat{\varphi}(x-y,\omega)\rho(y)dy,\quad x\in\mathbb{R}^d,\omega\in\Omega,
\end{equation}that satisfies 
$\tilde{\varphi}\in\mathrm{Lip}^3_m(\bb{R}^d\times\Omega)$ with $\partial_i\tilde{\varphi}\in 
\mathrm{Lip}^2_m(\bb{R}^d\times\Omega),i=1,\ldots d,$ where $m>0$ depends on the mollifier. A possible candidate for $\varphi$ is then $\tilde{\varphi}/m.$

\section\appendixname\subsection
{Invariant ergodic measure}\label{d}

The aim of the Appendix is to show the existence of an invariant ergodic measure for the process of the 
environment viewed from the particle defined in (\ref{11}) for a special class of diffusions in random environment 
and derive a formula 
for the limiting velocity of the diffusion, see Theorem \ref{t01}. We consider a class of diffusions that are 
Brownian motions which are perturbed by environment dependent drifts of the from given in (\ref{01}). 
Note that we will not assume any finite range 
dependence condition for the environment in order to derive the following results.\\

Let us consider a diffusion in random environment $\omega\in \Omega$ which is generated by the operator 
$\cal{L}^{\omega}$ given in 
(\ref{5}) with a drift of the form\begin{equation}\label{01}
b(x,\omega)\df \frac{\nabla\phi(x,\omega)}{2\phi(x,\omega)}+\frac{\Gamma(x,\omega)}{\phi(x,\omega)},\qquad 
x\in\bb{R}^d,\ \omega\in\Omega,
\end{equation}where $\phi\in \mathrm{Lip}^2_m(\bb{R}^d\times\Omega)$ 
for some $m>0,$ see below (\ref{19}) for the definition of $\mathrm{Lip}^2_m(\bb{R}^d\times\Omega)$, such that \begin{equation}\label{020}
{m}^{-1}\leq\inf_{x\in\bb{R}^d,\omega\in\Omega}\phi(x,\omega)\leq
\sup_{x\in\bb{R}^d,\omega\in\Omega}\phi(x,\omega)\leq m.
\end{equation}Moreover, we assume that $\phi(0,\omega)$ is a probability density with respect to the static 
measure of the environment $\bb{P}$ and 
$\Gamma=(\Gamma_1,\ldots,\Gamma_d)^T:\bb{R}^d\times\Omega\longrightarrow \bb{R}^d$ is a measurable function such that 
$\Gamma_i\in {\mathrm{Lip}^0_l}(\bb{R}^d\times\Omega)$ for some $l\geq 0$ 
and $\Gamma_i(\cdot,\omega)\in C^1(\bb{R}^d),\ i=1,\ldots,d,$ with $\nabla\cdot \Gamma(x,\omega)=0$ for all $x\in\bb{R}^d,\ 
\omega\in\Omega.$ Note 
that similarly as in the proof of (\ref{27}) one can show that the drift $b$ 
given in (\ref{01}) satisfies the Lipschitz condition (\ref{2}) and its Euclidean norm is uniformly bounded.
 Thus the martingale problem attached to (\ref{5}) with the above drift 
 starting from $x\in\bb{R}^d$ at time 0 is well-posed. Let in the sequel  
$P_{x,\omega},\ x\in\bb{R}^d,\ \omega\in\Omega,$ denote its unique solution and recall that $(X_t)_{t\geq 0}$ stands  for the 
coordinate process on $C(\bb{R}_+,\bb{R}^d).$ We can 
associate the canonical process on $\Omega$ defined by the environment $\omega\in\Omega$ as seen by an observer 
sitting on the particle, i.e.\begin{equation}\label{11}\left\{\begin{array}{rcl}
 \omega(t)&=&\tau_{X_t}(\omega),\\
 \omega(0)&=&\omega\in\Omega.\\
 \end{array}\right.  \end{equation}
 This map induces a measure $Q_{\omega}$ on the space of trajectories in $\Omega$ starting from $\omega\in\Omega$ 
 which can be shown to be a Markov process on $\Omega$ with semigroup given by$$
\cal{Q}_t \tilde{f}(\omega)\df E_{0,\omega}\left[\tilde{f}(\tau_{X_t}(\omega))\right],\quad t\geq 0, 
$$for all bounded measurable functions $\tilde{f}$ on $\Omega.$ Now we state the main theorem of the Appendix.

{\thm{\label{t01} The probability measure $d\bb{Q}\df \phi(0,\omega)d\bb{P}$ is an invariant ergodic measure for the Markov family 
$\{Q_{\omega}\}_{\omega\in\Omega}$ which is equivalent to 
the static measure $\bb{P}$ and 
\begin{equation}\label{10}
\bb{P}\times P_{0,\omega}= P_0\mbox{-a.s.,}\quad \lim_{t\to\infty}\frac{X_t}{t}=v \df E^{\bb{Q}}\left[b(0,\omega)\right].
\end{equation}
}}
\hspace*{-3pt}{\textbf{Proof:}} Since for all $x\in\bb{R}^d,\ \omega\in\Omega,\ 
\frac{1}{2}\Delta\phi(x,\omega)=\nabla\cdot\left(b(x,\omega)\phi(x,\omega)\right)
$ holds true, recall that $\nabla\cdot\Gamma(x,\omega)=0,$ and $\phi(0,\omega)$ 
is a probability density with respect to $\bb{P},$ we know that $\bb{Q}$ is an invariant 
probability measure for the family $\{Q_{\omega}\}_{\omega\in\Omega},$ 
see for instance Section 6 in \cite{va} and Section 6 in \cite{koz}. Due to the uniform ellipticity of
the diffusion matrix, which is in fact the identity matrix (see (\ref{5})),
it can be shown that $\bb{Q}$ is also ergodic. Indeed, a necessary and sufficient condition for ergodicity is that 
whenever a function $\tilde{g}$ on $\Omega$ which is square integrable with respect to the measure $\bb{Q}$ 
satisfies\begin{equation}\label{0100}
\cal{Q}_t\tilde{g}=\tilde{g},
\end{equation}$\bb{Q}$-a.s. for all $t>0,$ then $\tilde{g}$ has to be constant $\bb{Q}$-a.s., see 
Theorem 3.2.4 in \cite{da}. To show this we multiply both sides of (\ref{0100}) by $\tilde{g}$ and average 
over the environment with respect to the measure $\bb{Q}.$ After some manipulations using 
the invariance property of $\bb{Q}$ 
we find that (\ref{0100}) leads to\begin{equation}\label{0101}
E^{\bb{Q}}\left[E_{0,\omega}\left[\left(\tilde{g}(\tau_{X_t}(\omega))-\tilde{g}(\omega)\right)^2
\right]\right]=
E^{\bb{P}}\left[\phi(0,\omega)\int_{\bb{R}^d}p_{\omega}(t,0,y)\left(\tilde{g}(\tau_{y}(\omega))-\tilde{g}(\omega)
\right)^2dy\right]=0.
\end{equation}Due to the fact that $p_{\omega}(t,0,y)>0,$ for all $\omega\in\Omega, y\in\bb{R}^d,t>0,$ 
see Theorem 1 on page 67 in \cite{oleinik}, and $\bb{P}[\phi(0,\omega)\geq m^{-1}]=1,$ which comes from the uniform lower bound given in (\ref{020}), we can deduce from (\ref{0101}) that
\begin{equation}\label{0102}\bb{P}\mbox{-a.s.,}\qquad
\tilde{g}(\tau_{y}(\omega))=\tilde{g}(\omega)\qquad\mbox{for a.e. }y\in\bb{R}^d.
\end{equation}An application of the spatial ergodic theorem shown in \cite{dunschwartz}, see Theorem 10 on page 694, and the 
stationarity of the environment then show that (\ref{0102}) holds true for all $y\in\bb{R}^d.$
From the assumed ergodicity of the family of transformations $\{\tau_{x}\,:\,x\in\bb{R}^d\},$ mentioned above (\ref{1}),
it follows that $\tilde{g}$ is constant $\bb{P}$-a.s. and since $\bb{Q}$ is a measure equivalent to $\bb{P}$ 
due to (\ref{020}), $\tilde{g}$ is also constant $\bb{Q}$-almost surely. This shows the ergodicity of $\bb{Q}.$
Therefore, for 
each $\omega\in\Omega$ we have that for all bounded measurable function $\tilde{f}$ on $\Omega,$
\begin{equation}\label{0103}
\bb{Q}\times Q_{\omega}\mbox{-a.s.,}\qquad\lim_{t\to\infty}\frac{1}{t}\int_{0}^t\tilde{f}(\omega(s))ds=E^{\bb{Q}}[\tilde{f}],
\end{equation}which follows from Theorem 3.3.1 in \cite{da}. By definition of 
$Q_{\omega}$ we have that under the measure $Q_{\omega}$ the process $(\tilde{f}(\omega(s)))_{s\geq 0}$ has the 
same law as the process $(f(X_s,\omega))_{s\geq 0}$ under the measure $P_{0,\omega},$ where we define 
$f(x,\omega)\df\tilde{f}(\tau_x(\omega))$ for all $x\in\bb{R}^d,\ \omega\in\Omega,$ according to the definition (\ref{1}). Thus, 
(\ref{0103}) is equivalent to \begin{equation}\label{0104}
 \bb{Q}\times P_{0,\omega}\mbox{-a.s.,}\qquad
\lim_{t\to\infty}\frac{1}{t}\int_0^t f(X_s,\omega)ds\overset{(\ref{0104})}{=}E^{\bb{Q}}\left[
f(0,\omega)\right].
\end{equation}Since for $\omega\in\Omega,\ P_{0,\omega}$-a.s., $
X_t=\int_0^t b(X_s,\omega)ds+W_t,$ for all $t\geq 0,$
for some Brownian motion $(W_t)_{t\geq 0},$ (\ref{0104}) together with the law of large numbers for Brownian 
motion, see page 104 in \cite{ks}, yield
\begin{equation}\label{041}
\bb{Q}\times P_{0,\omega}\mbox{-a.s.,}\qquad
\lim_{t\to\infty}\frac{X_t}{t}=\lim_{t\to\infty}\frac{1}{t}\int_0^t b(X_s,\omega)ds=E^{\bb{Q}}\left[
b(0,\omega)\right].
\end{equation}Since $\bb{Q}$ and $\bb{P}$ are equivalent, (\ref{041}) holds true
$\bb{P}\times P_{0,\omega}$-almost surely. This concludes the proof of Theorem \ref{t01}.
\begin{flushright}$\Box$\end{flushright}

{\textbf{Acknowledgements:}} Let me thank my advisor Prof. A.-S. Sznitman
 for his constant advice during the completion of this work. I am also grateful for many helpful
 discussions with L. Goergen and D. Windisch.

\end{document}